\documentclass[review]{elsarticle}

\usepackage{lineno,hyperref}

\usepackage{amsbsy, amssymb, amsmath}
\usepackage{amsthm}

\modulolinenumbers[5]

\journal{https://arxiv.org/}

\usepackage{numcompress}\biboptions{authoryear}

\newtheorem{theorem}{Theorem} 
\newtheorem{lemma}{Lemma}[section]
\newtheorem{corollary}{Corollary} 
\newtheorem{example}{Example} 
\newtheorem{remark}{Remark} 

\def \be {\begin{eqnarray}}
\def \ee {\end{eqnarray}}

\def \bol{\boldsymbol{l}}
\def \boa {\boldsymbol{a}}

\def \bxm { \begin{example} }
\def \exm { \end{example} }
\def \brem { \begin{remark} }
\def \erem { \end{remark} }

\def\blem { \begin{lemma} }
\def\elem { {$\square$} \end{lemma} }
\def\btm { \begin{theorem} }
\def\etm { {$\blacksquare$} \end{theorem} }
\def\bpr { \begin{proof} }
\def\epr { \end{proof} }
\def\bcor { \begin{corollary} }
\def\ecor { {$\square$} \end{corollary} }

\begin{document}

\begin{frontmatter}

\title{\bf{A Proof of the Explicit Formula for\\Product Moments of\\
Multivariate Gaussian Random Variables}}

\author{Iickho Song\corref{mycorrespondingauthor}}
\cortext[mycorrespondingauthor]{Corresponding author. Tel.: +82-42-350-3445. This
study was supported by the National Research Foundation 
of Korea, under Grant NRF-2015R1A2A1A01005868, for which the author wishes to express his thanks.}
\ead{i.song@ieee.org}




\address{School of Electrical Engineering\\
Korea Advanced Institute of Science and Technology\\
291 Daehag Ro, Yuseong Gu,
Daejeon 34141 Korea}

\begin{abstract}
A detailed proof of a recent
result on explicit formulae for the product moments
$E \left \{ X_1^{a_1} X_2^{a_2} \cdots X_n^{a_n}\right \}$
of multivariate Gaussian random
variables is provided in this note.
\end{abstract}

\begin{keyword}
jointly Gaussian random variables
\sep Gaussian random vector
\sep joint moment 
\sep multivariate Gaussian random variables
\sep product moment
\end{keyword}

\end{frontmatter}

\linenumbers

\section{Introduction}
Recently, explicit formulae for the product moments of
multivariate Gaussian random variables have been derived in \cite{spl-is}. The formulae
not only are simpler and more concise than any of the other ones
but also allow us to promptly compute any term of the product moments.

Since the publication of the results, a huge number of
requests has been received for the details of the proof.
In this note, we provide a detailed proof of the main result of \cite{spl-is}.

\section{Proof of the Explicit Formula}

\subsection{The Explicit Formula}

Let $\left \{ X_i \right \}_{i=1}^{n}$
be multivariate Gaussian random
variables with mean vector
$\boldsymbol{\mu} = \left ( \mu_1, \mu_2, \cdots, \mu_n \right )$
and $n  \times  n$
variance-covariance matrix
$\Phi = \left[ \varphi_{ij} \right ]$,
where
$\mu_i = E \left \{ X_i \right \}$ denotes the expected value of $X_i$ and
$\varphi_{ij} = E \left \{ X_i X_j \right \}
- \mu_i \mu_j$ denotes the covariance
between $X_i$ and $X_j$ if $i \ne j$ and variance of $X_i$ if $i=j$.

One of the main results of \cite{spl-is} is as follows:

%
%
%
\btm
For the product moments
%
%
\be
E \left \{ \prod\limits_{i=1}^{n} X_i^{a_i} \right \} =
%
\sum\limits_{\bol \in S_{\boa} }
%
%
d_{\boa,\bol}
\left (
\prod\limits_{i=1}^{n}
\prod\limits_{j=i}^{n}
\varphi_{ij}^{l_{ij}} \right )
\left ( \prod\limits_{j=1}^{n}
\mu_{j}^{L_{\boa,j}} \right )
\label{eq-gen-99-1}
\ee
of multivariate Gaussian random variables, we have
%
%
%
%
%
%
%
%
%
%
\be
d_{\boa,\bol} =
\frac{\prod\limits_{k=1}^{n} a_k !}
{2^{M_{\bol}} \left ( \prod\limits_{i=1}^{n} \prod\limits_{j=i}^{n} l_{ij} ! \right )
\left ( \prod\limits_{j=1}^{n} L_{\boa,j} ! \right ) }
\label{eq-gen-99-2}
\ee
for $\bol \in S_{\boa}$, where $M_{\bol} = \sum\limits_{k=1}^{n} l_{kk}$.
\label{th-gen-99}
\etm

Here, the elements of the two sets
$\boa = \left \{ a_1, a_2, \cdots ,a_n \right \}$ and
$\bol= \left \{ \left \{ l_{ij} \right \}_{i=1}^{n}\right \}_{j=i}^{n}$
%
%
are 
nonnegative integers,
%
%
\be
L_{\boa,k} = a_k - l_{kk} - \sum\limits_{j=1}^{n}l_{jk}
\ee
%
with $l_{ji} = l_{ij}$, 
and 
$S_{\boa}$ denotes the collection 
of
%
%
all the sets $\bol$ such that
\be
%
%
%
\left \{ L_{\boa,k} \ge 0  \right \}_{k=1}^{n}.
%
%
\label{eq-abc-1-c-1}
\ee

\subsection{Proof}


Let us prove Theorem \ref{th-gen-99} by mathematical induction.
First, it is clear from Lemma 2.2, 
Theorem 1, and Theorem 2 of \cite{spl-is} 
that (\ref{eq-gen-99-2}) holds for $n=1, 2,$ and $3$, respectively.

Next, assume that (\ref{eq-gen-99-2}) 
is true when $n=m-1$. Then,
since
\be
E \left\{ X_1^{a_1} X_2^{a_2} \cdots X_{m}^{a_{m}}\right\}
= E \left\{ X_1^{a_1} X_2^{a_2} \cdots X_{m-1}^{a_{m-1}}\right\} E \left\{ X_{m}^{a_{m}}\right\}
\ee
when $\varphi_{1m} = \varphi_{2m} = \cdots =\varphi_{m-1,m} =0$,
we get
$\left[d_{\boldsymbol a_2, \boldsymbol l_2}\right]_{l_{1m}\mapsto 0, l_{2m}\mapsto 0, \cdots, l_{m-1,m} \mapsto 0} = d_{\boldsymbol a_1, \boldsymbol l_1} d_{a_m, l_{mm}}$,
or
\be
&& \hspace{-1.5cm}
\left[d_{\boldsymbol a_2, \boldsymbol l_2}\right]_{l_{1m}\mapsto 0, l_{2m}\mapsto 0, \cdots, l_{m-1,m} \mapsto 0}
= \nonumber \\
&&\frac{\prod\limits_{k=1}^{m} a_k!}
{2^{M_{\boldsymbol l_2}} \left(\prod\limits_{i=1}^{m-1} \prod\limits_{j=i}^{m-1} l_{ij}!\right) l_{mm}!
\left(\prod\limits_{j=1}^{m-1} L_{\boldsymbol a_1,j}!\right) \left(a_m-2l_{mm}\right)!}
\label{thm3-init}
\ee
using (\ref{eq-gen-99-2}) with $n=m-1$ and 
\be
d_{a_1 , l} = \frac{a_1!}{2^{l}l!\left(a_1-2l\right)!}
\label{eq-c2n-2}
\ee
(Eq. (5), \cite{spl-is}) for $l = 0, 1, \cdots, \lfloor \frac{a_1}{2} \rfloor $,
%
where
$\boldsymbol a_1= \left\{a_i\right\}_{i=1}^{m-1}$,
$\boldsymbol a_2= \boldsymbol a_1 \cup \left\{ a_{m} \right\}$,
$\boldsymbol l_1= \left\{ l_{ij}\right\}$ with $j=i,i+1,\cdots,m-1$ for $i=1,2,\cdots,m-1$, and
$\boldsymbol l_2= \boldsymbol l_1 \cup \left\{ l_{im}\right\}_{i=1}^{m}$.
Here, the symbol $\mapsto$ denotes a substitution:
For example, $\alpha \mapsto \beta$ means the substitution of
$\alpha$ with $\beta$.

Next, employing 
the result
%
\be
\frac {\partial}{\partial\varphi_{ij}}
E \left \{ \prod\limits_{k=1}^{n} X_k^{a_k} \right \} =
a_i a_j E \left \{ X_i^{a_i-1} X_j^{a_j-1} \underset{k \ne i, j}
{\prod\limits_{k=1}^{n}} X_k^{a_k} \right \}
\label{eq-price-92}
\ee
%
for $(i,j)=(1,m),(2,m),\cdots,(m-1,m)$ in \cite{price}, and then taking some steps
(similar to those from (16) to (20) in \cite{spl-is}), 
we will get
\be
\left[d_{\boldsymbol a_2, \boldsymbol l_2}\right]_{l_{im}\mapsto l_{im}+1}
=\frac{a_i!a_m!\left[d_{\boldsymbol a_2, \boldsymbol l_2}\right]_{l_{im}\mapsto 0, a_i \mapsto a_i-l_{im}-1, a_m \mapsto a_m-l_{im}-1}}
{\left(a_i-l_{im}-1\right)!\left(a_m-l_{im}-1\right)!\left(l_{im}+1\right)!}
\label{thm3-l-im-plus}
\ee
and then, by substituting $l_{im}+1$ with $l_{im}$ in (\ref{thm3-l-im-plus}),
we have
\be
d_{\boldsymbol a_2, \boldsymbol l_2}
=\frac{a_i!a_m!\left[d_{\boldsymbol a_2, \boldsymbol l_2}\right]_{l_{im}\mapsto 0, a_i \mapsto a_i-l_{im}, a_m \mapsto a_m-l_{im}}}
{\left(a_i-l_{im}\right)!\left(a_m-l_{im}\right)!\left(l_{im}\right)!}
\label{thm3-coeff-0}
\ee
for $i=1,2,\cdots,m-1$.
Now, letting $l_{2m}=0$ in (\ref{thm3-coeff-0}) with $i=1$, we get
\be
\left[d_{\boldsymbol a_2, \boldsymbol l_2}\right]_{l_{2m}\mapsto 0}
=\frac{a_1!a_m!\left[d_{\boldsymbol a_2, \boldsymbol l_2}\right]_{l_{1m}\mapsto 0, l_{2m}\mapsto 0, a_1 \mapsto a_1-l_{1m}, a_m \mapsto a_m-l_{1m}}}
{\left(a_1-l_{1m}\right)!\left(a_m-l_{1m}\right)! l_{1m}!},
\label{thm3-l-2m-zero}
\ee
which can be used into 
(\ref{thm3-coeff-0}) with $i=2$ to produce 
\be
&& \hspace{-1cm}
d_{\boldsymbol a_2, \boldsymbol l_2}
=\frac{a_2!a_m!}{\left(a_2-l_{2m}\right)!\left(a_m-l_{2m}\right)!l_{2m}!}  \times \frac{a_1! \left(a_m-l_{2m}\right)!}{\left(a_1-l_{1m}\right)!\left(a_m-l_{1m}-l_{2m}\right)! l_{1m}!} \nonumber \\
&& \quad  \times \left[d_{\boldsymbol a_2, \boldsymbol
l_2}\right]_{l_{1m}\mapsto 0, l_{2m}\mapsto 0, a_1 \mapsto
a_1-l_{1m}, a_2 \mapsto a_2-l_{2m}, a_m \mapsto
a_m-l_{1m}-l_{2m}}. \label{thm3-coeff-1} \ee Subsequently, letting
$l_{3m}=0$ in (\ref{thm3-coeff-1}), we get \be && \hspace{-1cm}
\left[d_{\boldsymbol a_2, \boldsymbol l_2}\right]_{l_{3m}\mapsto
0} = \frac{a_1!a_2!a_m!}
{l_{1m}!l_{2m}! \left(a_1-l_{1m}\right)!\left(a_2-l_{2m}\right)! \left(a_m-l_{1m}-l_{2m}\right)!} \nonumber \\
&&  \times \left[d_{\boldsymbol a_2, \boldsymbol
l_2}\right]_{l_{km}\mapsto 0\text{ for }k=1,2,3,
 \text{ }a_1 \mapsto a_1-l_{1m}, a_2 \mapsto a_2-l_{2m}, a_m \mapsto a_m-l_{1m}-l_{2m}},
\label{thm3-l-3m-zero}
\ee
which can be employed 
into (\ref{thm3-coeff-0}) with $i=3$ to produce 
\be
d_{\boldsymbol a_2, \boldsymbol l_2}
&=&\frac{a_3!a_m!}{\left(a_3-l_{3m}\right)!\left(a_m-l_{3m}\right)!l_{3m}!} \nonumber \\
&&  \times \frac{a_1!a_2! \left(a_m-l_{3m}\right)!}
{l_{1m}!l_{2m}! \left(a_1-l_{1m}\right)!\left(a_2-l_{2m}\right)! \left(a_m-l_{1m}-l_{2m}-l_{3m}\right)!} \nonumber \\
&&  \times \left[d_{\boldsymbol a_2, \boldsymbol
l_2}\right]_{l_{km}\mapsto 0, a_k \mapsto a_k-l_{km}\text{ for
}k=1,2,3,\text{ }a_m \mapsto a_m-l_{1m}-l_{2m}-l_{3m}}.
\label{thm3-coeff-2} \ee
From
(\ref{thm3-coeff-2}), it can be formulated \be d_{\boldsymbol a_2,
\boldsymbol l_2} &=& \frac{\left(\prod\limits_{k=1}^i a_k! \right)
a_m!} {\left(\prod\limits_{k=1}^{i} l_{km}!\right)
 \left\{ \prod\limits_{k=1}^{i} \left(a_k-l_{km}\right)! \right\} \left(a_m-
 \sum\limits_{k=1}^i l_{km}\right)!} \nonumber \\
&&  \times \left[d_{\boldsymbol a_2, \boldsymbol
l_2}\right]_{l_{km}\mapsto 0, a_k \mapsto a_k-l_{km}\text{ for
}k=1,2,3,\text{ }a_m \mapsto a_m-\sum\limits_{k=1}^i l_{km}}
\label{thm3-coeff-gen} \ee with $i=3$. Next, letting $l_{4m}=0$ in
(\ref{thm3-coeff-gen}), we will get $\left[d_{\boldsymbol a_2,
\boldsymbol l_2}\right]_{l_{4m}\mapsto 0}$ which can be employed
into (\ref{thm3-coeff-0}) with $i=4$ to produce
(\ref{thm3-coeff-gen}) with $i=4$. After proceeding the above
procedures until $i=m-1$, we will eventually get \be
d_{\boldsymbol a_2, \boldsymbol l_2} &=&
\frac{\left(\prod\limits_{k=1}^{m-1} a_k!\right) a_m!}
{\left(\prod\limits_{k=1}^{m-1} l_{km}!\right) \left \{
\prod\limits_{k=1}^{m-1} \left(a_k-l_{km}\right)! \right\}
\left(a_m-
\sum\limits_{k=1}^{m-1} l_{km}\right)!} \nonumber \\
&&  \times \left[d_{\boldsymbol a_2, \boldsymbol
l_2}\right]_{l_{km}\mapsto 0, a_k \mapsto a_k-l_{km}\text{ for
}k=1,2,\cdots,m-1,\text{ }a_m
 \mapsto a_m-\sum\limits_{k=1}^{m-1} l_{km}}
\label{thm3-coeff-m-1}
\ee
from $\left[d_{\boldsymbol a_2, \boldsymbol l_2}\right]_{l_{m-1,m}\mapsto 0}$
(obtained by letting $l_{m-1,m}=0$ in (\ref{thm3-coeff-gen}) with $i=m-2$) and
(\ref{thm3-coeff-0}) with $i=m-1$.
Finally, by combining (\ref{thm3-init}) into (\ref{thm3-coeff-m-1}), we can get
\be
&& \hspace{-2cm}
d_{\boldsymbol a_2, \boldsymbol l_2} = \frac{\left(\prod\limits_{k=1}^{m-1} a_k!\right) a_m!}
{\left(\prod\limits_{k=1}^{m-1} l_{km}!\right)
\left\{ \prod\limits_{k=1}^{m-1} \left(a_k-l_{km}\right)! \right\}
\left(a_m- \sum\limits_{k=1}^{m-1} l_{km}\right)!} \nonumber \\
&& \hspace{-1cm}  \quad   \times \frac{\left\{
\prod\limits_{k=1}^{m-1} \left(a_k-l_{km}\right)! \right\}
\left(a_m- \sum\limits_{k=1}^{m-1} l_{km}\right)!}
{2^{M_{\boldsymbol l_2}} \left(\prod\limits_{k=1}^{m-1}
\prod\limits_{j=i}^{m-1} l_{kj}!\right) l_{mm}!
%
} \nonumber \\
&& \hspace{-1cm}  \quad   \times \frac{1}
{\prod\limits_{j=1}^{m-1} \left(\left[L_{\boldsymbol
a_1,j}\right]_{a_k \mapsto a_k-l_{km}\text{ for }k=1,2,\cdots,m-1}
\right)!} \nonumber \\
&& \hspace{-1cm}  \quad   \times \frac{1}{\left(a_m-\sum\limits_{k=1}^{m-1} l_{km}-2l_{mm}\right)!} \nonumber \\
&& \hspace{-1cm}
=
\frac{\prod\limits_{k=1}^{m} a_k!}
{2^{M_{\boldsymbol l_2}} \left(\prod\limits_{i=1}^{m} \prod\limits_{j=i}^{m} l_{ij}!\right)
\left(\prod\limits_{j=1}^{m} L_{\boldsymbol a_2,j}!\right)},
\ee
which implies that (\ref{eq-gen-99-2}) 
holds when $n=m$ also.


\section{Conclusion}

We have provided a detailed proof of an elegant
explicit formulae for product moments
$E \left \{ X_1^{a_1} X_2^{a_2}
\cdots X_n^{a_n} \right \}$
of multivariate Gaussian variables derived recently.



\end{document}